\documentclass[10pt,journal]{IEEEtran}

\usepackage{amssymb}
\usepackage[cmex10]{amsmath}
\usepackage{multirow}
\usepackage{graphicx}
\usepackage{epstopdf}
\usepackage{amsthm}

\graphicspath{{./fig/}}

\newcommand{\ve}{\mathbf}
\newcommand{\m}{\mathbf}

\newcommand{\vef}[1]{\mathbf{\tilde{\mathbf{#1}}}} 

\newtheorem{Proposition}{Proposition}

\title{CWCU LMMSE Estimation: Prerequisites and Properties}
\author{Mario Huemer, ~\emph{Senior Member,~IEEE}, Oliver Lang, \emph{Student Member,~IEEE}  \\
\authorblockA{ Johannes Kepler University Linz, Institute of Signal Processing \\ 
							 Altenbergerstr. 69, 4040 Linz, Austria \{mario.huemer, oliver.lang\}@jku.at}}

\begin{document}
%

\maketitle
%
\begin{abstract}
The classical unbiasedness condition utilized e.g. by the best linear unbiased estimator (BLUE) is very stringent. By softening the ''global'' unbiasedness condition and introducing component-wise conditional unbiasedness conditions instead, the number of constraints limiting the estimator's performance can in many cases significantly be reduced. In this work we investigate the component-wise conditionally unbiased linear minimum mean square error (CWCU LMMSE) estimator for different model assumptions. The prerequisites in general differ from the ones of the LMMSE estimator. We first derive the CWCU LMMSE estimator under the jointly Gaussian assumption of the measurements and the parameters. Then we focus on the linear model and discuss the CWCU LMMSE estimator for jointly Gaussian parameters, and for mutually independent (and otherwise arbitrarily distributed) parameters, respectively. In all these cases the CWCU LMMSE estimator incorporates the prior mean and the prior covariance matrix of the parameter vector. For the remaining cases optimum linear CWCU estimators exist, but they may correspond to globally unbiased estimators that do not make use of prior statistical knowledge about the parameters. Finally, the beneficial properties of the CWCU LMMSE estimator are demonstrated with the help of a well-known channel estimation application. 
\end{abstract}

\begin{keywords}
	Estimation, Bayesian Estimation, Best Linear Unbiased Estimator, BLUE, Linear Minimum Mean Square Error, LMMSE, CWCU, Channel   
	Estimation.
\end{keywords}

\section{Introduction}
\label{sec:Introduction}
\let\thefootnote\relax\footnotetext{This work was supported by the Austrian Science Fund (FWF): I683-N13.}

Usually, when we talk about unbiased estimation of a parameter vector $\ve{x}\in\mathbb{C}^{n\times 1}$ out of a measurement vector $\ve{y}\in\mathbb{C}^{m\times 1}$, then the estimation problem is treated in the classical framework \cite{Kay-Est.Theory}, \cite{Mueller_2012}. Letting $\hat{\ve{x}} = \ve{g}(\ve{y})$ be an estimator of $\ve{x}$, then the classical unbiased constraint asserts that 
\begin{equation}
	E_\ve{y}[\hat{\ve{x}}] = \int{\ve{g}(\ve{y}) p(\ve{y};\ve{x}) d\ve{y} = \ve{x}} \hspace{0.4cm} \text{for all possible } \ve{x}, \label{equ:cuLMMSE045}
\end{equation}
where $p(\ve{y};\ve{x})$ is the probability density function (PDF) of vector $\ve{y}$ parametrized by the unknown parameter vector $\ve{x}$. The index of the expectation operator shall indicate the PDF over which the averaging is performed. 
Eq.~\eqref{equ:cuLMMSE045} can also be formulated in the Bayesian framework, where the parameter vector $\ve{x}$ is treated as random. Here, the corresponding problem arises by demanding global conditional unbiasedness, i.e. 
\begin{equation}
	E_{\ve{y}|\ve{x}}[\hat{\ve{x}}|\ve{x}] = \int{\ve{g}(\ve{y}) p(\ve{y}|\ve{x}) d\ve{y}} = \ve{x} 
									\hspace{0.4cm} \text{for all possible } \ve{x}.\label{equ:cuLMMSE046}
\end{equation}
The attribute \emph{global} indicates that the condition is made on the whole parameter vector 
$\ve{x}$. However, the constricting requirement in \eqref{equ:cuLMMSE046} prevents the exploitation of prior knowledge about the parameters, and hence leads to a significant reduction in the benefits brought about by the Bayesian framework. 

In component-wise conditionally unbiased (CWCU) Bayesian parameter estimation \cite{Slock-Oct2005}-\cite{Huemer_2014_Asilomar}, instead of constraining the estimator to be globally unbiased, we aim for achieving conditional unbiasedness on one parameter component at a time. Let $x_i$ be the $i^{th}$ element of $\ve{x}$, 
and $\hat{x}_i = g_i(\ve{y})$ be the estimator of $x_i$. Then the CWCU constraints are
\begin{equation}
	E_{\ve{y}|x_i}[\hat{x}_i|x_i] = 
									\int{g_i(\ve{y}) p(\ve{y}|x_i)} d\ve{y} = x_i, \label{equ:cuLMMSE047}  
\end{equation}
for all possible $x_i$ ($i=1,2,...,n$). 
The CWCU constraints are less stringent than the global conditional unbiasedness condition in \eqref{equ:cuLMMSE046}, and it will turn out that a CWCU estimator in many cases allows the incorporation of prior knowledge about the statistical properties of the parameter vector.

The paper is organized as follows: In Section \ref{sec:GaussianAssumption} we discuss the prerequisites and the solution of the CWCU linear minimum mean square error (LMMSE) estimator under the jointly Gaussian assumption of $\ve{x}$ and $\ve{y}$. 
In Section \ref{sec:LinearModel} we discuss the CWCU LMMSE estimator under the linear model assumption, and we extend the findings of \cite{Slock-Oct2005}. Here we particularly distinguish between jointly Gaussian, and mutually independent (and otherwise arbitrarily distributed) parameters. Finally, in Section \ref{sec:IEEE_CIR_Model} the CWCU LMMSE estimator is compared against the best linear unbiased estimator (BLUE) and the LMMSE estimator in a well-known channel estimation example.

\section{CWCU LMMSE Estimation under the jointly Gaussian Assumption} \label{sec:GaussianAssumption}
We assume that a vector parameter $\ve{x}\in\mathbb{C}^{n\times 1}$ is to be estimated based on a measurement vector 
$\ve{y}\in\mathbb{C}^{m\times 1}$. As in LMMSE estimation we constrain the estimator to be linear (or actually affine), such that 
\begin{equation} 
	\hat{\ve{x}} = \m{E} \ve{y} + \ve{c}, \hspace{0.5cm} \m{E}\in\mathbb{C}^{n\times m},  \ve{c}\in\mathbb{C}^{n\times 1}.  \label{equ:cuLMMSE049e}
\end{equation} 
Note that in LMMSE estimation no assumptions on the specific form of the joint PDF $p(\ve{x},\ve{y})$ have to be made. However, the situation is different in CWCU LMMSE estimation. Let us consider the $i^{th}$ component of the estimator 
\begin{equation}
	\hat{x}_i = \ve{e}_i^H \ve{y} +c_i, 
	\label{equ:cuLMMSE050b}
\end{equation} 
where $\ve{e}_i^H$ denotes the $i^{th}$ column of the estimator matrix $\m{E}$. The conditional mean of $\hat{x}_i$ can be written as
\begin{equation}
	E_{\ve{y}|x_i}[\hat{x}_i|x_i] = 
				\ve{e}_i^H E_{\ve{y}|x_i}[\ve{y}|x_i] + c_i. \label{equ:cuLMMSE050a}
\end{equation} 
A closer inspection of \eqref{equ:cuLMMSE050a} reveals that $E_{\ve{y}|x_i}[\hat{x}_i|x_i] = x_i$ can be fulfilled for all possible $x_i$ if the conditional mean 
$E_{\ve{y}|x_i}[\ve{y}|x_i]$ is a linear function of $x_i$. For jointly Gaussian $\ve{x}$ and $\ve{y}$ this is the case and we have
\begin{equation}
	E_{\ve{y}|x_i}[\ve{y}|x_i] = E_{\ve{y}}[\ve{y}] + 
									(\sigma_{x_i}^2)^{-1}\m{C}_{\ve{y} x_i} (x_i -E_{x_i}[x_i]),
\end{equation}
where $\m{C}_{\ve{y} x_i} = E_{\ve{y}, x_i}[(\ve{y} - E_{\ve{y}}[\ve{y}])(x_i - E_{x_i}[x_i])^H]$, and $\sigma_{x_i}^2$ is the variance of $x_i$. $E_{\ve{y}|x_i}[\hat{x}_i|x_i] = x_i$ is fulfilled if
\begin{equation}
	\ve{e}_i^H \m{C}_{\ve{y} x_i} = \sigma_{x_i}^2   \label{equ:CWCULMMSE010}
\end{equation}
and 
\begin{equation}
	c_i = E_{x_i}[x_i] - \ve{e}_i^H E_{\ve{y}}[\ve{y}]. \label{equ:CWCULMMSE011}
\end{equation}
Inserting \eqref{equ:cuLMMSE050b}, \eqref{equ:CWCULMMSE010} and \eqref{equ:CWCULMMSE011} in the Bayesian MSE cost function $E_{\ve{y},\ve{x}}[|\hat{x}_i - x_i|^2]$ immediately leads to the optimization problem
\begin{equation}
	\ve{e}_{\text{CL},i} = \mathrm{arg}\mathop{\mathrm{min}}_{\ve{e}_i} \hspace{0.1cm} 
															    \left(\ve{e}_i^H \m{C}_{\ve{y}\ve{y}} \ve{e}_i - \sigma_{x_i}^2 \right)
																	\hspace{0.25cm} \text{s.t.}  \hspace{0.1cm}	\ve{e}_i^H \m{C}_{\ve{y} x_i}  = \sigma_{x_i}^2,
\end{equation}
where ''CL'' shall stand for CWCU LMMSE. The solution can be found with the Lagrange multiplier method and is given by
\begin{equation}
	\ve{e}_{\text{CL},i} = \frac{\sigma_{x_i}^2}{\m{C}_{x_i\ve{y}}\m{C}_{\ve{y}\ve{y}}^{-1}\m{C}_{\ve{y}x_i}}
												 \m{C}_{\ve{y}\ve{y}}^{-1}\m{C}_{\ve{y}x_i}. \label{equ:CWCULMMSE012}
\end{equation}
Using $\m{E}_{\text{CL}} = [e_{\text{CL},1}, e_{\text{CL},2}, \hdots, e_{\text{CL},n}]^H$ together with \eqref{equ:CWCULMMSE011} and \eqref{equ:CWCULMMSE012} immediately leads us to the first part of the 
\smallskip

\begin{Proposition} \label{prop:CWCULMMSE001} If $\ve{x}\in\mathbb{C}^{n\times 1}$ and $\ve{y}\in\mathbb{C}^{m\times 1}$ are jointly Gaussian then the CWCU LMMSE estimator minimizing the Bayesian MSEs $E_{\ve{y},\ve{x}}[|\hat{x}_i - x_i|^2]$ under the constraints 
$E_{\ve{y}|x_i}[\hat{x}_i|x_i]=x_i$ for $i = 1,2,\cdots,n$
is given by 
\begin{equation}
	\hat{\ve{x}}_{\mathrm{CL}} = E_{\ve{x}}[\ve{x}] +\m{E}_{\mathrm{CL}} (\ve{y} - E_{\ve{y}}[\ve{y}]), \label{equ:CWCULMMSE001}
\end{equation}
with 
\begin{equation}
	\m{E}_{\mathrm{CL}} = \m{D}\m{C}_{\ve{x}\ve{y}}\m{C}_{\ve{y}\ve{y}}^{-1}, \label{equ:CWCULMMSE002} 
\end{equation}
where the elements of the real diagonal matrix $\m{D}$ are 
\begin{equation}
	[\m{D}]_{i,i} = \frac{\sigma_{x_i}^2}{\m{C}_{x_i\ve{y}}\m{C}_{\ve{y}\ve{y}}^{-1}\m{C}_{\ve{y}x_i}}. \label{equ:CWCULMMSE003}
\end{equation}
The mean of the error $\ve{e} = \ve{x}-\hat{\ve{x}}_{\mathrm{CL}}$ (in the Bayesian sense) is zero, and the error covariance matrix 
$\m{C}_{\ve{e}\ve{e},\mathrm{CL}}$ which is also the minimum Bayesian MSE matrix $\m{M}_{\hat{\ve{x}}_{\mathrm{CL}}}$ is 
\begin{equation}
	\m{C}_{\ve{e}\ve{e},\mathrm{CL}} = \m{M}_{\hat{\ve{x}}_{\mathrm{CL}}} = 
				\m{C}_{\ve{x}\ve{x}} - \m{A}\m{D} - \m{D}\m{A} + \m{D}\m{A}\m{D},
\end{equation}
with $\m{A} = \m{C}_{\ve{x}\ve{y}}\m{C}_{\ve{y}\ve{y}}^{-1}\m{C}_{\ve{y}\ve{x}}$. The minimum Bayesian MSEs are 
$\mathrm{Bmse}(\hat{x}_{\mathrm{CL},i}) = [\m{M}_{\hat{\ve{x}}_{\mathrm{CL}}}]_{i,i}$.
\end{Proposition}
\smallskip

The part on the error performance can simply be proofed by inserting in the definition of $\ve{e}$ and $\m{C}_{\ve{e}\ve{e}}$, respectively. From \eqref{equ:CWCULMMSE002} it can be seen that the CWCU LMMSE estimator matrix can be derived as the product of the diagonal matrix $\m{D}$ with the LMMSE estimator matrix $\m{E}_{\text{L}}=\m{C}_{\ve{x}\ve{y}}\m{C}_{\ve{y}\ve{y}}^{-1}$. Furthermore, we have 
$E_{\ve{y}|x_i}[\hat{x}_{\mathrm{L},i}|x_i] = [\m{D}]_{i,i}^{-1} x_i + (1-[\m{D}]_{i,i}^{-1})E_{x_i}[x_i]$ for the LMMSE estimator. It can be shown that $\m{D}$ can also be written as
\begin{equation}
	\m{D} = \text{diag}\{ \m{C}_{\ve{x}\ve{x}}\} \left(  \text{diag}\{ \m{A} \} \right)^{-1}. \label{equ:CWCULMMSE005a}
\end{equation} 

The CWCU LMMSE estimator will in general not commute over linear transformations, an exception is discussed in \cite{Huemer_2014_Asilomar}.

\section{CWCU LMMSE Estimation under Linear Model Assumptions} \label{sec:LinearModel}

In the following it will be seen that some of the prerequisites of Proposition \ref{prop:CWCULMMSE001} can be relaxed when incorporating details of the data model into the derivation of the estimator. From now on we limit our considerations to the linear model 
\begin{equation} 
	\ve{y} = \m{H}\ve{x} + \ve{n}, \label{equ:cuLMMSE049b}
\end{equation}
where $\m{H}\in\mathbb{C}^{m\times n}$ is a known observation matrix, $\ve{x}\in \mathbb{C}^{n\times 1}$ is a parameter vector with mean $E_{\ve{x}}[\ve{x}]$ and covariance matrix $\m{C}_{\ve{x}\ve{x}}$, and $\ve{n}\in \mathbb{C}^{m\times 1}$ is a zero mean noise vector with covariance matrix $\m{C}_{\ve{n}\ve{n}}$ and independent of $\ve{x}$. Additional assumptions on 
$\ve{x}$ and $\ve{n}$ will vary in the following. We note that the CWCU LMMSE estimator for the linear model under the assumption of white Gaussian noise has already been derived in \cite{Slock-Oct2005}. 

\subsection{Solution for correlated Gaussian parameters} \label{sec:LinearModelGaussian}

For the linear model the covariance matrices required in \eqref{equ:CWCULMMSE002} and \eqref{equ:CWCULMMSE003} become $\m{C}_{\ve{y}\ve{y}} = \m{H}\m{C}_{\ve{x}\ve{x}}\m{H}^H + \m{C}_{\ve{n}\ve{n}}$, $\m{C}_{\ve{x}\ve{y}} = \m{C}_{\ve{x}\ve{x}}\m{H}^H$, $\m{C}_{x_i\ve{y}} = \m{C}_{x_i\ve{x}}\m{H}^H$ and $\m{C}_{\ve{y}x_i} = \m{H}\m{C}_{\ve{x}x_i}$.
If the assumptions made on the linear model above hold and if $\ve{x}$ and $\ve{n}$ are both Gaussian, then they are jointly Gaussian. Furthermore, since $[\ve{x}^T,\ve{y}^T]^T$ is a linear transformation of $[\ve{x}^T,\ve{n}^T]^T$, $\ve{x}$ and $\ve{y}$ are jointly Gaussian, too. We could therefore simply insert 
the above covariance matrices
into the equations given in Proposition \ref{prop:CWCULMMSE001}. However, the jointly Gaussian assumption for $\ve{x}$ and $\ve{n}$ can significantly be relaxed. This can be shown by incorporating the linear model assumption already earlier in the derivation of the estimator. Let
$\ve{h}_i\in\mathbb{C}^{m\times 1}$ be the $i^{th}$ column of $\m{H}$, $\bar{\m{H}}_i\in\mathbb{C}^{m\times (n-1)}$ the matrix resulting from $\m{H}$ by deleting $\ve{h}_i$, and $\bar{\ve{x}}_i\in\mathbb{C}^{(n-1)\times 1}$ the vector resulting from $\ve{x}$ after deleting $x_i$. Then we can write the $i^{th}$ component of $\hat{\ve{x}}$ in the form
\begin{align}
	\hat{x}_i &= \ve{e}_i^H \ve{y} +c_i = \ve{e}_i^H (\ve{h}_i x_i + \bar{\m{H}}_i \bar{\ve{x}}_i + \ve{n}) + c_i. \label{equ:CWCULMMSE015} 
\end{align}
The conditional mean of $\hat{x}_i$ becomes
\begin{equation}
	E_{\ve{y}|x_i}[\hat{x}_i|x_i] = \ve{e}_i^H\ve{h}_i x_i + 
					\ve{e}_i^H\bar{\m{H}}_i E_{\bar{\ve{x}}_i|x_i}[\bar{\ve{x}}_i|x_i] + c_i. \label{equ:CWCULMMSE016}
\end{equation} 
The CWCU constraint $E_{\ve{y}|x_i}[\hat{x}_i|x_i] =x_i$ for a particular $i$ can be fulfilled if
$E_{\bar{\ve{x}}_i|x_i}[\bar{\ve{x}}_i|x_i]$ is a linear function of $x_i$. This is true if $\ve{x}$ is Gaussian, i.e.
$\ve{x} \sim \mathcal{CN}(E_{\ve{x}}[\ve{x}], \m{C}_{\ve{x}\ve{x}})$, since then we have 
$E_{\bar{\ve{x}}_i|x_i}[\bar{\ve{x}}_i|x_i]=E_{\bar{\ve{x}}_i}[\bar{\ve{x}}_i]+
(\sigma_{x_i}^2)^{-1}\m{C}_{\bar{\ve{x}}_i x_i}(x_i-E_{x_i}[x_i])$. 
Note that the only requirement on the noise vector so far was its independence on $\ve{x}$. 
Following similar arguments as above we end up at the constrained optimization problem 
\begin{align}
	\ve{e}_{\text{CL},i} &= \mathrm{arg}\mathop{\mathrm{min}}_{\ve{e}_i} \hspace{0.1cm} 
													\left(\ve{e}_i^H (\m{H}\m{C}_{\ve{x}\ve{x}}\m{H}^H + \m{C}_{\ve{n}\ve{n}}) \ve{e}_i - 
											    \sigma_{x_i}^2 \right) \nonumber \\
											 &	\hspace{1.9cm} \text{s.t.}  \hspace{0.1cm}	\ve{e}_i^H \m{H}\m{C}_{\ve{x}x_i} = \sigma_{x_i}^2.
													\label{equ:CWCULMMSE017}
\end{align}
Solving \eqref{equ:CWCULMMSE017} leads to 

\begin{Proposition} \label{prop:CWCULMMSE003} If the observed data $\ve{y}$ follow the linear model in \eqref{equ:cuLMMSE049b}, where $\ve{y}\in\mathbb{C}^{m\times 1}$ is the data vector, $\m{H}\in\mathbb{C}^{m\times n}$ is a known observation matrix, $\ve{x}\in\mathbb{C}^{n\times 1}$ is a parameter vector with prior PDF $\mathcal{CN}(E_{\ve{x}}[\ve{x}], \m{C}_{\ve{x}\ve{x}})$, and $\ve{n}\in\mathbb{C}^{m\times 1}$ is a zero mean noise vector with covariance matrix $\m{C}_{\ve{n}\ve{n}}$ and independent of $\ve{x}$ (the joint PDF of $\ve{x}$ and $\ve{n}$ is otherwise arbitrary), then the CWCU LMMSE estimator minimizing the Bayesian MSEs $E_{\ve{y},\ve{x}}[|\hat{x}_i - x_i|^2]$ under the constraints 
$E_{\ve{y}|x_i}[\hat{x}_i|x_i]=x_i$ for $i = 1,2,\cdots,n$ is given by \eqref{equ:CWCULMMSE001} with
\begin{equation}
	\m{E}_{\mathrm{CL}} = \m{D}\m{C}_{\ve{x}\ve{x}}\m{H}^H(\m{H}\m{C}_{\ve{x}\ve{x}}\m{H}^H+\m{C}_{\ve{n}\ve{n}})^{-1}, 
											\label{equ:CWCULMMSE004}
\end{equation}
where the elements of the real diagonal matrix $\m{D}$ are 
\begin{equation}
	[\m{D}]_{i,i} = \frac{\sigma_{x_i}^2}{\m{C}_{x_i\ve{x}}\m{H}^H(\m{H}\m{C}_{\ve{x}\ve{x}}\m{H}^H+\m{C}_{\ve{n}\ve{n}})^{-1}
									\m{H}\m{C}_{\ve{x}x_i}}. \label{equ:CWCULMMSE005}
\end{equation} 
\end{Proposition}
\smallskip

The same formulas would result by inserting 
for the covariance matrices in the equations given by Proposition \ref{prop:CWCULMMSE001}, however, the prerequisites in Proposition 1 and 2 differ. Also the error measures can formally be derived by inserting in the equations of Proposition \ref{prop:CWCULMMSE001}. 

\subsection{Solution for mutually independent parameters}

For mutually independent parameters it is possible to further relax the prerequisites on $\ve{x}$. In this situation \eqref{equ:CWCULMMSE016} becomes  
\begin{equation}
	E_{\ve{y}|x_i}[\hat{x}_i|x_i] = \ve{e}_i^H\ve{h}_i x_i + 
					\ve{e}_i^H\bar{\m{H}}_i E_{\bar{\ve{x}}_i}[\bar{\ve{x}}_i] + c_i, 
\end{equation} 
since $E_{\bar{\ve{x}}_i|x_i}[\bar{\ve{x}}_i|x_i]$ is no longer a function of $x_i$. The CWCU constraints are fulfilled if 
$\ve{e}_i^H\ve{h}_i=1$ and $c_i = -\ve{e}_i^H\bar{\m{H}}_i E_{\bar{\ve{x}}_i}[\bar{\ve{x}}_i]$, and no further assumptions on the PDF of $\ve{x}$ are required. Again following similar arguments as above we end up at a constrained optimization problem \cite{Huemer_2014_Asilomar}. Solving it leads to

\begin{Proposition} \label{prop:CWCULMMSE004} If the observed data $\ve{y}$ follow the linear model where $\ve{y}\in\mathbb{C}^{m\times 1}$ is the data vector, $\m{H}\in\mathbb{C}^{m\times n}$ is a known observation matrix, $\ve{x}\in\mathbb{C}^{n\times 1}$ is a parameter vector with mean $E_{\ve{x}}[\ve{x}]$, mutually independent elements and covariance matrix $\m{C}_{\ve{x}\ve{x}} = \mathrm{diag}\{\sigma_{x_1}^2,\sigma_{x_2}^2,\cdots,\sigma_{x_n}^2\}$, $\ve{n}\in\mathbb{C}^{m\times 1}$ is a zero mean noise vector with covariance matrix $\m{C}_{\ve{n}\ve{n}}$ and independent of $\ve{x}$ (the joint PDF of $\ve{x}$ and $\ve{n}$ is otherwise arbitrary), then the CWCU LMMSE estimator minimizing the Bayesian MSEs $E_{\ve{y},\ve{x}}[|\hat{x}_i - x_i|^2]$ under the constraints 
$E_{\ve{y}|x_i}[\hat{x}_i|x_i]=x_i$ for $i = 1,2,\cdots,n$ is given by \eqref{equ:CWCULMMSE001} and \eqref{equ:CWCULMMSE004}, where the elements of the real diagonal matrix $\m{D}$ are 
\begin{equation}
	[\m{D}]_{i,i} = \frac{1}{\sigma_{x_i}^2\ve{h}_i^H(\m{H}\m{C}_{\ve{x}\ve{x}}\m{H}^H+\m{C}_{\ve{n}\ve{n}})^{-1}\ve{h}_i}.
										\label{equ:CWCULMMSE006}
\end{equation}
\end{Proposition}
\smallskip

In \cite{Huemer_2014_Asilomar} we showed that for mutually independent parameters $\ve{e}_{\text{CL},i}$ is independent of $\sigma_{x_i}^2$ and also given by
$\ve{e}_{\text{CL},i} = (\ve{h}_i^H \m{C}_i^{-1} \ve{h}_i)^{-1} \m{C}_i^{-1} \ve{h}_i$,
where $\m{C}_i= \bar{\m{H}}_i \m{C}_{\bar{\ve{x}}_i\bar{\ve{x}}_i} \bar{\m{H}}_i^H + \m{C}_{\ve{n}\ve{n}}$. 
Furthermore, we showed that  
\begin{equation}
	[\m{D}]_{i,i} = (\ve{e}_{\text{L},i}^H \ve{h}_i)^{-1}, \label{equ:CWCU_Journal034}
\end{equation}
where $\ve{e}_{\text{L},i}^H$ is the $i^{th}$ row of the LMMSE estimator. It therefore holds that 
$\mathrm{diag}\{\m{E}_{\text{CL}}\m{H}\} = \ve{1}$.


\subsection{Other cases}

If $\ve{x}$ is whether Gaussian nor a parameter vector with mutually independent parameters, then we have the following possibilities: If $E_{\ve{y}|x_i}[\hat{x}_i|x_i]$ is a linear function of $x_i$ for all $i = 1,2,\cdots,n$ then we can derive the CWCU LMMSE estimator similar as in Section \ref{sec:LinearModelGaussian}. In the remaining cases still an estimator can be found that fulfills the CWCU constraints. By studying \eqref{equ:CWCULMMSE016} it can be seen that the choice $\ve{e}_i^H\ve{h}_i = 1$, $\ve{e}_i^H\bar{\m{H}}_i=\ve{0}$ together with
$c_i=0$ for all $i = 1,2,\cdots,n$ will do the job. Inserting these constraints into the Bayesian MSE cost functions and solving the constrained optimization problem leads to 
\begin{equation}
	\hat{\ve{x}}_{\text{B1}} = \m{E}_{\text{B1}}\ve{y} = 
				(\m{H}^H\m{C}_{\ve{n}\ve{n}}^{-1}\m{H})^{-1} \m{H}^H\m{C}_{\ve{n}\ve{n}}^{-1} \ve{y}, \label{equ:CWCULMMSE013}
\end{equation}
with $\m{C}_{\ve{e}\ve{e},\text{B1}} = (\m{H}^H\m{C}_{\ve{n}\ve{n}}^{-1}\m{H})^{-1}$ as the Bayesian error covariance matrix.
$\ve{e}_i^H\ve{h}_i = 1 \hspace{0.1cm} \& \hspace{0.1cm} \ve{e}_i^H\bar{\m{H}}_i=\ve{0}$ is equivalent to $\m{E}\m{H}=\m{I}$. This implies $\hat{\ve{x}}_{\text{B1}} = \m{E}_{\text{B1}}\ve{y}= \ve{x} + \m{E}_{\text{B1}}\ve{n}$. It is seen that the estimator in 
\eqref{equ:CWCULMMSE013} fulfills the global unbiasedness condition 
$E_{\ve{y}|\ve{x}}[\hat{\ve{x}}_{\text{B1}}|\ve{x}]=\ve{x}$ for every $\ve{x}\in\mathbb{C}^{n\times 1}$. This estimator, which is the BLUE if $\ve{x}$ can be any vector in $\in\mathbb{C}^{n\times 1}$, does not account for any prior knowledge about $\ve{x}$. In some situations a better globally unbiased estimator exists. 
If for example it is known that 
$\ve{x}$ lies in a linear subspace of $\mathbb{C}^{n\times 1}$ spanned by the columns of the full rank matrix $\m{V}\in\mathbb{C}^{n\times p}$ with $p<n$, such that $\ve{x}=\m{V}\ve{z}$, then each estimator with $\m{E}\m{H}\m{V}=\m{V}$ fulfills 
$E_{\ve{y}|\ve{x}}[\hat{\ve{x}}|\ve{x}]=\ve{x}$ for every $\ve{x}=\m{V}\ve{z}$ lying in the subspace \cite{Huemer12_1}. However, 
$\m{E}\m{H}\m{V}=\m{V}$ is less stringent than $\m{E}\m{H}=\m{I}$, so inserting this weaker constraints into the Bayesian MSE cost functions and solving the constrained optimization problem leads to a better performing estimator which is
\begin{equation}
	\hat{\ve{x}}_{\text{B2}} = \m{E}_{\text{B2}}\ve{y} = 
			\m{V}(\m{V}^H\m{H}^H\m{C}_{\ve{n}\ve{n}}^{-1}\m{H}\m{V})^{-1} \m{V}^H\m{H}^H\m{C}_{\ve{n}\ve{n}}^{-1}\ve{y},
			\label{equ:CWCULMMSE014}
\end{equation}
with the Bayesian error covariance matrix
\begin{equation}
	\m{C}_{\ve{e}\ve{e},\text{B2}} = \m{V}(\m{V}^H\m{H}^H\m{C}_{\ve{n}\ve{n}}^{-1}\m{H}\m{V})^{-1} \m{V}^H. \label{equ:CWCULMMSE019}
\end{equation}
In case the assumption that $\ve{x}$ lies in the linear subspace defined by the full rank matrix $\m{V}$ holds, this estimator is in fact the BLUE. 

Note that the estimator in \eqref{equ:CWCULMMSE013} 
(and the one in \eqref{equ:CWCULMMSE014} once its underlying assumptions are fulfilled) 
can of course also be applied if the prerequisites in Proposition \ref{prop:CWCULMMSE003} are fulfilled. However, in this case the CWCU LMMSE estimator is always the one given in Proposition \ref{prop:CWCULMMSE003}.

\section{Example: Channel Estimation} \label{sec:IEEE_CIR_Model}

\begin{figure}[!t]
\centering
\includegraphics[width=3.0in]{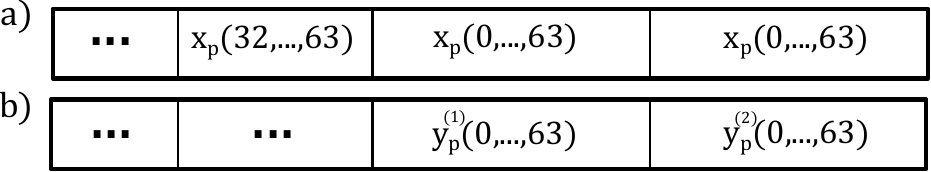}
\caption{a) Preamble including two long training symbols and a long guard interval for channel estimation; b) received long training symbols.}
\label{fig:preamble}
\end{figure}

As an application to demonstrate the properties of the CWCU LMMSE estimator we choose the well-known channel estimation problem for IEEE 802.11a/g/n WLAN standards \cite{IEEE99}. The standards define two identical preamble symbols $\ve{x}_p\in\mathbb{C}^{64\times 1}$, cf. Fig.~\ref{fig:preamble}a, which are designed such that the frequency domain version $\vef{x}_{p}=\m{F}_N\ve{x}_p$, shows $\pm 1$ at 52 subcarrier positions (indexes $\{1,...,26,38,...63\}$) and zeros at the remaining unused 12 subcarriers (indexes $\{0,27,...,37\}$). Here $\m{F}_N$ is the length $N=64$ discrete Fourier transform (DFT) matrix, and $\tilde{( \cdot)}$ denotes a vector in the frequency domain. 
With the carrier selection matrix 
$\m{B}\in\{0,1\}^{64 \times 52}$, cf. \cite{Huemer_2014_Asilomar},
the vector of nonzero (used) subcarrier symbols can be written as $\vef{x}_{p,u} = \m{B}^T \m{F}_N \ve{x}_p$. 
$\m{B}^T$ deletes the elements of $\vef{x}_{p}$ that correspond to the zero-subcarriers. We furthermore introduce the diagonal matrix $\m{D}_p = \mathrm{diag} \{ \vef{x}_{p,u} \}$ which fulfills $\m{D}_p^H \m{D}_p=\m{I}$ because of 
$\vef{x}_{p,u}\in\{-1,1\}^{52 \times 1}$.

The channel impulse response (CIR) is modeled as $\ve{h}\sim\mathcal{CN}(\ve{0}, \m{C}_{\ve{h}\ve{h}})$, with 
\begin{equation}
	\m{C}_{\ve{h}\ve{h}} = \mathrm{diag}\{\sigma_0^2,\sigma_1^2,...,\sigma_{l_h-1}^2\}, \label{equ:IEEE_CIR002}
\end{equation}
and exponentially decaying power delay profile according to 
$\sigma_i^2 = \left(1-\mathrm{exp}(-\tfrac{T_s}{\tau_{rms}})\right)\mathrm{exp}(-\tfrac{i T_s}{\tau_{rms}})$ for $i=0,1,...,l_h-1$.
Here, $l_h$ is the length of the CIR. $T_s$ and $\tau_{rms}$ are the sampling time and the channel delay spread, respectively, which are chosen as $T_s=50$ns, $\tau_{rms}=100$ns in our setup. Note that the channel length $l_h$ can be assumed to be considerably smaller than the DFT length $N$. In the following we assume $l_h=16$.

Let $\ve{y}_{p}^{(1)}$ and $\ve{y}_{p}^{(2)}$ be the two received, channel distorted time domain preamble symbols, cf. Fig.~\ref{fig:preamble}b, $\vef{y}_{p,u}^{(i)}= \m{B}^T \m{F}_N\ve{y}_{p}^{(i)}$ for $i=1,2$, and 
$\overline{\tilde{\ve{y}}}=\frac{1}{2}(\vef{y}_{p,u}^{(1)}+\vef{y}_{p,u}^{(2)})$. Then $\overline{\tilde{\ve{y}}}$ can be modeled as
\begin{align}
	\overline{\tilde{\ve{y}}} &= \m{D}_p \vef{h}_u + \vef{n} \label{equ:IEEE_CIR030} \\
														&= \m{D}_p \m{B}^T \vef{h} + \vef{n} \label{equ:CWCULMMSE025} \\
														&= \m{D}_p \m{B}^T \m{M}_1 \ve{h} + \vef{n}. \label{equ:IEEE_CIR032} 
\end{align}
Here $\vef{h}_u\in\mathbb{C}^{52\times 1}$ is the frequency response at the occupied subcarriers, $\vef{h}\in\mathbb{C}^{64\times 1}$ is the full length frequency response including the unused frequency bins, and $\vef{n}$ is a complex white Gaussian noise vector with covariance matrix $\m{C}_{\tilde{\ve{n}}\tilde{\ve{n}}}=(N\sigma_n^2/2) \m{I}$, where $\sigma_n^2$ is the time domain noise variance. $\m{M}_1\in\mathbb{C}^{64\times 16}$ consists of the first $l_h$ columns of $\m{F}_N$.

\begin{figure}[!t]
\centering
\includegraphics[width=3.3in]{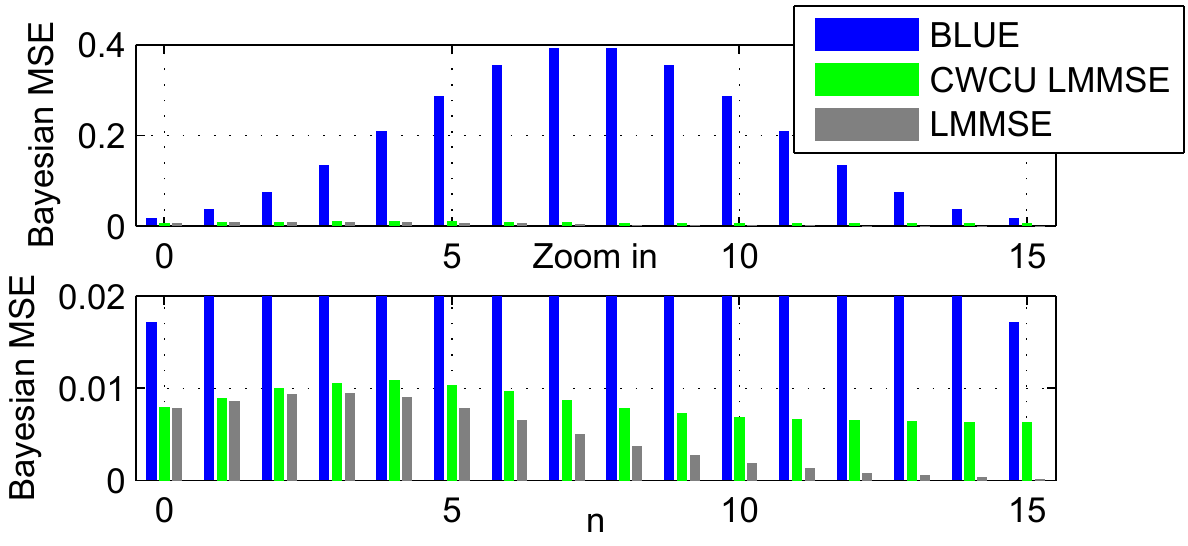}
\caption{Bayesian MSEs of the estimated CIR coefficients.}
\label{fig:CIR_time_domain}
\end{figure}

From \eqref{equ:IEEE_CIR032} the BLUE, the LMMSE and the CWCU LMMSE estimator for the channel impulse response $\ve{h}$ follow to
\begin{align}
	&\hat{\ve{h}}_{\text{B}} = \left(\m{M}_1^H \m{B} \m{B}^T \m{M}_1\right)^{-1} \m{M}_1^H \m{B}\m{D}_p^{-1}\overline{\tilde{\ve{y}}} \label{equ:IEEE_CIR015}	 \\
	&\hat{\ve{h}}_{\text{L}} = \left(\m{M}_1^H \m{B} \m{B}^T \m{M}_1 + \frac{N\sigma_n^2}{2} \m{C}_{\ve{h}\ve{h}}^{-1}\right)^{-1} \m{M}_1^H \m{B}\m{D}_p^{-1} \overline{\tilde{\ve{y}}}  \label{equ:IEEE_CIR016} \\
	&\hat{\ve{h}}_{\text{CL}} = \m{D} \hat{\ve{h}}_{\text{L}}.  \label{equ:IEEE_CIR017}
\end{align}
Here $\m{D}$ shall be used from Proposition \ref{prop:CWCULMMSE004}, since the elements of $\ve{h}$ are mutually independent.

Fig.~\ref{fig:CIR_time_domain} shows the Bayesian MSEs of the estimated CIR coefficients for the different estimators (for 
$\sigma_n^2=0.01$). The performance drawback of the BLUE mainly comes from the fact that no measurements are available at the large gap from subcarrier 27 to 37. The LMMSE estimator and the CWCU LMMSE estimator incorporate the prior knowledge from \eqref{equ:IEEE_CIR002} which results in a huge performance gain over the BLUE. The CWCU LMMSE estimator almost reaches the LMMSE estimator's performance, and in contrast to the latter it additionally shows the beneficial property of conditional unbiasedness. \smallskip

\begin{figure}[!t]
\centering
\includegraphics[width=3.3in]{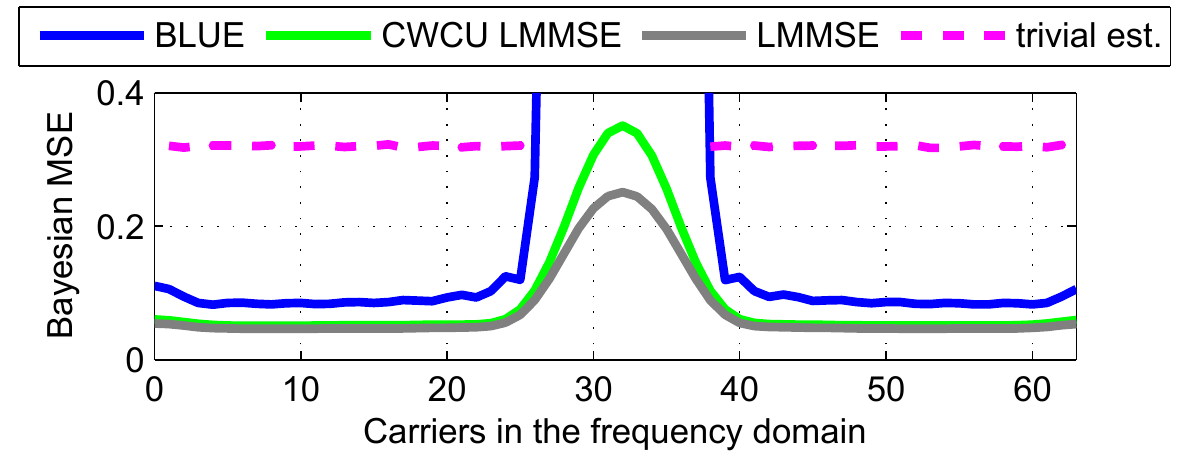}
\caption{Bayesian MSEs for the elements of $\hat{\vef{h}}_{\text{B}}$, $\hat{\vef{h}}_{\text{L}}$, 
$\hat{\vef{h}}_{\text{CL}}$, and $\hat{\vef{h}}_{u,\text{trivial}}$.}
\label{fig:complete_CIR_frequency_domain}
\end{figure}

We now turn to frequency response estimators. From \eqref{equ:IEEE_CIR030} a straight forward trivial estimator for $\vef{h}_u$ follows to $\hat{\vef{h}}_{u,\text{trivial}} = \m{D}_p^{-1}\overline{\tilde{\ve{y}}}$. This estimator fulfills the unbiasedness condition \eqref{equ:cuLMMSE045} for every $\vef{h}_u\in\mathbb{C}^{52\times 1}$, but since $\vef{h}_u$ lies in a linear subspace of $\mathbb{C}^{52\times 1}$ (spanned by the columns of $\m{B}^T \m{M}_1$), this estimator is not the BLUE. The LMMSE estimator which commutes over linear transformations is 
$\hat{\vef{h}}_{\text{L}} = \m{M}_1 \hat{\ve{h}}_{\text{L}}$. $\hat{\vef{h}}_{\text{B}} = \m{M}_1 \hat{\ve{h}}_{\text{B}}$ corresponds to the BLUE as discussed in \eqref{equ:CWCULMMSE014}. The CWCU LMMSE estimator does not commute over general linear transformations, but by using the prior covariance matrix $\m{C}_{\vef{h}\vef{h}}=\m{M}_1\m{C}_{\ve{h}\ve{h}}\m{M}_1^H$, $\hat{\vef{h}}_{\text{CL}}$ (which is not $\m{M}_1\hat{\ve{h}}_{\text{CL}}$) can easily be derived from $\hat{\vef{h}}_{\text{L}}$ by applying Proposition \ref{prop:CWCULMMSE003}. Fig.~\ref{fig:complete_CIR_frequency_domain} shows the Bayesian MSEs of $\hat{\vef{h}}_{\text{B}}$, $\hat{\vef{h}}_{\text{L}}$, $\hat{\vef{h}}_{\text{CL}}$, and $\hat{\vef{h}}_{u,\text{trivial}}$. Practically we are usually only interested in estimates at the occupied 52 subcarrier positions. However, in this work we study the estimators' performances at all 64 subcarrier positions, since this highlights some interesting properties. The BLUE significantly outperforms $\hat{\vef{h}}_{u,\text{trivial}}$, and in contrast to the latter it is able to estimate the frequency response at all subcarriers. However, the performance at the large gap from subcarrier 27 to 37, where no measurements are available, is extremely poor. (The maximum Bayesian MSE appears at subcarrier 32 and exhibits the huge value of around 36.) In contrast, the LMMSE estimator and the CWCU LMMSE estimator show excellent interpolation properties along the huge gap. As in the time domain the CWCU LMMSE estimator comes close to the LMMSE estimator's performance. 

\bigskip

\section{Conclusion}
In this work we investigated the CWCU LMMSE estimator under different model assumptions. First we derived the estimator for the case that the measurements and the parameters are jointly Gaussian. Then we concentrated on the linear model, where the only assumption made on the noise vector is its independence on the parameter vector. The CWCU LMMSE estimator has been derived for correlated Gaussian parameter vectors, and for the case that the parameters are mutually independent (and otherwise distributed arbitrarily). For the remaining cases the CWCU LMMSE estimator may correspond to a globally unbiased estimator.

\newpage



\end{document}